\documentclass[12pt,a4paper,reqno]{amsart}
\usepackage[utf8]{inputenc}
\usepackage[T1]{fontenc}
\usepackage{amsmath}
\usepackage{amsfonts}
\usepackage{amssymb}
\usepackage{mathtools}
\usepackage{graphicx}
\usepackage{xcolor}
\usepackage{comment}
\usepackage[left=3cm,right=3cm,bottom=2.8cm,top=2.8cm]{geometry}
\usepackage{calrsfs}
\usepackage{hyperref}
\hypersetup{
	colorlinks=true,
	linkcolor=blue,
	citecolor=blue,      
	urlcolor=black,
	linktocpage=true
}
\newtheorem{theorem}{Theorem}[section]
\newtheorem{proposition}[theorem]{Proposition}
\newtheorem{lemma}[theorem]{Lemma}
\newtheorem{corollary}[theorem]{Corollary}
\newtheorem{definition}[theorem]{Definition}

\newtheorem{remark}[theorem]{Remark}
\newtheorem{conjecture}[theorem]{Conjecture}
\numberwithin{equation}{section}

\DeclareMathAlphabet{\pazocal}{OMS}{zplm}{m}{n}

\def\N{\mathbb{N}}
\def\Z{\mathbb{Z}}
\def\Q{\mathbb{Q}}
\def\R{\mathbb{R}}
\def\C{\mathbb{C}}
\def\Cs{\mathbb{C}[\! [s] \!]}
\def\Ct{\mathbb{C}[\! [t] \!]}
\def\Kt{K[\! [t] \!]}
\def\Cst{\mathbb{C}[\![s,t]\!]}
\def\Kst{K[\![s,t]\!]}
\def\val{\mathrm{val}}
\def\ft{\widetilde{f}}
\def\A{\mathbb{A}}
\def\Spec{\operatorname{Spec}}
\def\codim{\operatorname{codim}}
\renewcommand{\O}{\mathcal{O}}
\def\I{\mathcal{I}}
\def\Xsm{X_\mathrm{sm}}
\def\Xs{X_\mathrm{sing}}
\def\CC{\mathcal{E}}
\renewcommand{\L}{\mathcal{L}}
\def\X{\mathcal{X}}
\def\J{\operatorname{Jac}}
\def\gt{\widetilde{\gamma}}
\def\F{\mathbb{F}}
\def\CC{\mathcal{C}}
\def\AC{\mathrm{AC}}
\def\NN{\mathcal{N}}
\def\lcm{\mathrm{lcm}}

\title{The embedded Nash problem in singular spaces: The case of surfaces}
\author{Javier de la Bodega}
\date{}

\begin{document}
	\maketitle
	
	\begin{abstract}
		We introduce the embedded Nash problem allowing singularities in the ambient space, and solve it in the case of surfaces, generalizing \cite[Theorem 1.22]{BdlB}. 
	\end{abstract}
	
	\section{Introduction}
	
	Let $X$ be a complex algebraic variety, $Z$ a closed subvariety of $X$, and $\Sigma$ a Zariski-closed subset of $Z$. For each integer $m \geq 1$, the $m$-contact locus of $(X,Z,\Sigma)$ is the set $\X_m$ of arcs in $X$ based in $\Sigma$ that have intersection multiplicity $m$ with $Z$. These subsets are central in motivic integration, in the definition of the motivic zeta function, and in the monodromy conjecture \cite{DL98, DL99}. When $X$ is smooth, and $Z$ is a hypersurface with an isolated singularity, contact loci are also conjecturally related to the Floer homology of iterates of the monodromy as predicted by the arc-Floer conjecture \cite{BBLN,BL}.
	
	Under the assumption of $X$ being smooth, contact loci were studied in \cite{ELM}. Moreover, in Remark 2.8 of loc. cit. the problem of determining its irreducible components in terms of an embedded resolution was posed. This is known as the \emph{embedded Nash problem}, and it was addressed in depth in \cite{BdlB}. In particular, the problem was completely solved in the case of $X$ being smooth and $Z$ being an unibranch curve. This was a key step when proving the arc-Floer conjecture for plane curves in \cite{BL}.
	
	In any case, there is no reason why we should stick to the assumption of $X$ being smooth, and in fact, the embedded Nash problem could be stated in greater generality. In this paper, we present the embedded Nash problem for 
	singular ambient spaces, and solve it when $X$ is a surface, $Z$ is an effective Cartier divisor and $\Sigma$ is a closed point in $Z$. For a precise statement, see Theorem \ref{main}.
	
	The proof adapts the ideas of \cite{FdB} to the embedded setting, combining algebro-geometric and topological techniques. In this sense, the proof reminds of the solution of the classical Nash problem for surfaces by Fernández de Bobadilla and Pe Pereira \cite{BP}. De Fernex and Docampo later solved the problem by purely algebro-geometric arguments that involved MMP constructions \cite{FD}. An analogous proof for the embedded Nash problem for 
	surfaces should also be expected.
	
	This work generalizes \cite[Theorem 1.22]{BdlB}. However, in loc. cit. the obtained result is stronger, for when $X$ is a smooth surface the irreducible components of the contact loci are disjoint. We expect this result to remain valid under possibly mild assumptions on $X$, but the techniques used in the present work are not enough to clarify this point. We leave it as an open problem in Conjecture \ref{conj}.

	Hopefully, as in the smooth ambient case, solving the embedded Nash problem for surfaces will guide us toward an extension of the arc-Floer conjecture where 
	singularities are allowed in the ambient space.	
	\\\\
	\noindent
	\textbf{Acknowledgements.} The author is grateful to Javier Fernández de Bobadilla for suggesting the problem and for the useful conversations they had during the development of the paper. The author is also grateful to Nero Budur for carefully reading the first versions of the paper and suggesting changes and improvements. Finally, María Pe Pereira observed that a minimality condition was missing in a previous version of the paper, for which the author is especially thankful.
	
	The author was supported by the grants PRE2019-087976, PID2020-114750GB-C33, SEV-2023-2026 from the Spanish Ministry of Science and G097819N, G0B3123N from FWO.
	
	\section{embedded Nash problem}
	
	We first recall some basic notions about arc spaces.
	
	\subsection{Arc spaces}
	Let $X$ be a $\C$-scheme. Given a field extension $K$ of $\C$, a $\C$-morphism $\Spec \Kt \to X$ is called a \textbf{$K$-arc} in $X$. The functor that associates to a $\C$-scheme $S$ the set $X(S \times \Ct)$ is representable, and the $\C$-scheme representing it will be denoted by $\L(X)$. It is known as the \textbf{arc space} of $X$. It is generally not locally of finite type, but it is locally of countable type if so is $X$. This construction is functorial,  i.e. every $\C$-morphism $h: Y \to X$ of $\C$-schemes induces a $\C$-morphism $h_\infty: \L(Y) \to \L(X)$.
	
	By definition, for every field extension $K$ of $\C$, the set of $K$-valued points of $\L(X)$ is the set of $K$-arcs in $X$. A point in $\L(X)$ will be called a \textbf{schematic arc} in $X$. 
	A $K$-arc corresponding to a schematic arc $\gamma \in \L(X)$ will be called a \textbf{representative arc} of $\gamma$.
	
	A \textbf{$K$-wedge} is a $K$-arc in the space of arcs of $X$,  i.e. a $\C$-morphism $\Spec \Kst \to X$. Analogously, the space of arcs of the space of arcs of $X$ is known as the \textbf{space of wedges} of $X$, and will be denoted by $\L^{(2)}(X)$. In particular, the set of its $K$-valued points is $X(\Kst)$.
	
	Let $\gamma$ be a schematic arc in $X$. Pick a representative arc $\Spec \Kt \to X$ of $\gamma$. The image of the closed point of $\Spec \Kt$ will be called the \textbf{base point} of $\gamma$ and will be denoted by $\gamma(0)$. Similarly, the image of the generic point of $\Spec \Kt$ will be called the \textbf{generic point} of $\gamma$, and will be denoted by $\gamma(\eta)$. These definitions do not depend on the choice of the representative arc of $\gamma$, and they define continuous maps $\L(X) \to X$. When $X$ is irreducible, we say that a schematic arc is \textbf{fat} if its generic point is the generic point of $X$.
	
	Let $\gamma$ be a schematic arc in $X$ and $Z$ a closed subscheme of $X$. Pick a representative arc $\Spec \Kt \to X$ of $\gamma$. The scheme-theoretic inverse of $Z$ is a subscheme of $\Spec \Kt$ defined by the ideal $\langle t^m \rangle$ for some $m \in \N \cup \{\infty\}$ (here, we set $t^\infty = 0$). The integer $m$ is known as the \textbf{intersection multiplicity} of $\gamma$ and $Z$, and will be denoted by $\gamma \cdot Z$. This definition does not depend on the choice of representative arc of $\gamma$. The intersection multiplicity is positive if and only if the base point of $\gamma$ is contained in $Z$, and it is finite if and only if the generic point of $\gamma$ is not in $Z$.
	
	A $\C$-arc $\Spec \Ct \to X$ is said to be \textbf{convergent} if it factors through the canonical morphism $\Spec \Ct \to \Spec \C \{t\}$, where $\C \{t\}$ denotes the ring of convergent power series. Similarly, we define the notion of a \textbf{convergent} $\C$-wedge.
	
	Recall that a power series $\gamma(t) \in \Ct$ is said to be \textbf{algebraic} if it is the root in $x$ of a nonzero polynomial $P(t,x) \in \C[t,x]$. The set of algebraic power series is a subring of $\Ct$, and we denote it by $A$. A $\C$-arc $\Spec \Ct \to X$ is said to be \textbf{algebraic} if it factors through the canonical morphism $\Spec \Ct \to \Spec A$. Similarly, we define the notion of a \textbf{algebraic} $\C$-wedge. We remark that algebraic $\C$-arcs (resp. algebraic $\C$-wedges) are convergent.
	
	The arc space of $X$ is typically not of finite type, not even  when $X$ is. Instead, it can be approximated by schemes of finite type. More precisely, for every $\ell \in \N$, the functor that associates to a $\C$-scheme $S$ the set $X(S \times \C[t]/t^{\ell+1})$ is representable, and the $\C$-scheme representing it will be denoted by $\L_\ell(X)$. It is known as the \textbf{$\ell$-jet space} of $X$, and it is of finite type if so is $X$.
	
	Again the jet scheme construction is functorial,  i.e. every $\C$-morphism $h: Y \to X$ of $\C$-schemes induces a $\C$-morphism $h_\ell: \L_\ell(Y) \to \L_\ell(X)$
	
	For every $\ell \geq \ell'$ there are natural $\C$-morphisms $\tau^\ell_{\ell'}: \L_\ell(X) \to \L_{\ell'}(X)$ and $\tau_\ell: \L(X) \to \L_\ell(X)$. The induced $\C$-morphism
	\[
	\lim_{\longleftarrow} \L_\ell(X) \longrightarrow \L(X)
	\]
	is an isomorphism by \cite[Theorem 1.1]{Bh}.
	
	Suppose that $X$ is of finite type, and let $\ell \in \N$. We say that a subset $C$ of $\L(X)$ is a \textbf{cylinder of level $\ell$} if there exists a constructible subset $C_\ell$ of $\L_\ell(X)$ such that $C = \tau_\ell^{-1}(C_\ell)$. We say $C$ is a \textbf{cylinder} if it is a cylinder of level $\ell$ for some $\ell \in \N$. Note that if $C$ is a cylinder of level $\ell'$, then it is also a cylinder of level $\ell$ for every $\ell \geq \ell'$.
	
	\subsection{Contact loci} Let $X$ be an integral separated $\C$-scheme of finite type of dimension $d$, $Z$ a closed subscheme of $X$ and $\Sigma$ a nonempty closed subset of $Z$. We will always assume that $X$ is smooth away from $Z$.
	\begin{definition}
		For each positive integer $m \geq 1$, we call the subset
		\[
		\X_m \coloneqq \X_m(X,Z,\Sigma) \coloneqq \{\gamma \in \L(X) \ | \ \gamma \cdot Z = m, \ \gamma(0) \in \Sigma\}
		\]
		the \textbf{$m$-contact locus} of $(X,Z,\Sigma)$.
	\end{definition}
	When $\Sigma = Z$, we will only write $\X_m(X,Z)$ for $\X_m(X,Z,Z)$. We will also denote
	\[
	\X_{\geq m} \coloneqq \X_{\geq m}(X,Z,\Sigma) \coloneqq \{\gamma \in \L(X) \ | \ \gamma \cdot Z \geq m, \ \gamma(0) \in \Sigma\}.
	\]
	The subset $\X_{\geq m}$ is a closed cylinder of $\L(X)$, and since $\X_m = \X_{\geq m} - \X_{\geq m+1}$, we conclude that $\X_m$ is a locally closed cylinder of $\L(X)$.
	
	Note that contact loci are local objects, in the sense that if $U$ is an open subset of $X$ containing $\Sigma$, then
	\[
	\X_m(X,Z,\Sigma) = \X_m(U,U \cap Z,\Sigma).
	\]
	Therefore, when $\Sigma$ is a point, we may assume that $X$ is affine.
	
	Let $h: Y \to X$ be an embedded resolution of $(X,Z,\Sigma)$, i.e. a proper birational morphism such that $Y$ is smooth over $\C$, the map $h$ is an isomorphism over $X-Z$, and $h^{-1}(Z)$ and $h^{-1}(\Sigma)$ are simple normal crossing divisors. Note that if we had not assumed that $X$ is smooth away from $Z$, then there would not exist any embedded resolution.
	
	Let $\CC$ be the set of irreducible components of $h^{-1}(Z)$, and $\CC_\Sigma$ the set of irreducible components of $h^{-1}(\Sigma)$. Since $\Sigma \subseteq Z$, every component of $h^{-1}(\Sigma)$ is also a component of $h^{-1}(Z)$, i.e. $\CC_\Sigma \subseteq \CC$.
	
	We write the scheme-theoretic inverse image of $Z$ as
	\[
	h^{-1}(Z) = \sum_{E \in \CC} N_EE.
	\]
	
	For each divisor $E \in \CC$, we write
	\[
	E^\circ \coloneqq E - \bigcup (\CC-\{E\}),
	\]
	and denote by $\NN_E$ the set of divisors of $\CC-\{E\}$ meeting $E$.

	By the valuative criterion of properness, every schematic arc $\gamma \in \X_m$ admits a unique lift $\gt$, i.e. a schematic arc $\gt \in \L(Y)$ such that $h_\infty(\gt)=\gamma$. Since $\gamma(0) \in \Sigma$, there exists a divisor $E \in \CC_\Sigma$ containing $\gt(0)$. Thus, for every $E \in \CC_\Sigma$, we consider the subset
	\[
	\X_{m,E} \coloneqq \{\gamma \in \X_m \ | \ \gt(0) \in E\}.
	\]
	The situation simplifies substantially if we impose an extra condition on $h$.
	
	\begin{definition}[\cite{McL}]
		We say $h$ is \textbf{$m$-separating} if given two different intersecting divisors $E,E' \in \CC$, we have that $N_E+N_{E'} > m$.
	\end{definition}
	
	Under this condition, the subsets $\X_{m,E}$ and $\X_{m,E'}$ are disjoint if $E, E'\in \CC_\Sigma$ are different divisors. Also, $\X_{m,E}$ is nonempty if and only if $N_E$ divides $m$. In that case, we will say that $E$ is an \textbf{$m$-divisor}. Overall, we have the union
	\begin{equation} \label{union}
		\X_m = \bigcup_{\text{$E$ $m$-divisor}} \X_{m,E} \ , 
	\end{equation}
	which is set-theoretically disjoint.
	
	\subsection{Connection to maximal divisorial sets}
	The subsets $\X_{m,E}$ are closely related to maximal divisorial sets, see \cite[Definition 2.8]{Is08}. Indeed, every divisorial valuation $v$ on $X$ admits a maximal closed irreducible subset $C_X(v)$ of $\L(X)$ inducing the valuation $v$.
	
	For every $m$-divisor $E$, denote by $\val_E: K(X)^\times \to \Z$ the divisorial valuation on $X$ induced by $E$, and set $w_E \coloneqq \frac{m}{N_E}\val_E$. The connection is the following.
	
	\begin{proposition} \label{generic-point}
		Let $E$ be an $m$-divisor. Then the generic point of $C_X(w_E)$ lies in $\X_{m,E}$. In particular, the closure of $\X_{m,E}$ in $\L(X)$ is $C_X(w_E)$.
	\end{proposition}
	\begin{proof}
		By \cite[Propostion 3.4]{Is08}, the generic point of $C_X(w_E)$ is of the form $\gamma = h_\infty(\gt)$, where $\gt \cdot E = \frac{m}{N_E}$. Hence $\gamma \cdot Z \geq m$. Since $E$ is an $m$-divisor, there exists $\delta \in \X_{m,E}$, so in particular, $\delta \in C_X(w_E)$. Since $\gamma$ is its generic point and the intersection multiplicity is upper-semicontinuous, we have $m = \delta \cdot Z \geq \gamma \cdot Z \geq m$, so $\gamma \cdot Z = m$ as desired.
	\end{proof}

	\subsection{The subsets $\X_{m,E}$ are cylinders}
	
	In \cite[Theorem A]{ELM}, it is proven that if $X$ is smooth, then $\X_{m,E}$ is a cylinder for every $m$-divisor $E$. We will show that the same holds when $X$ is allowed to have singularities, which will be an easy application of \cite{EM}. Before, we recall some notions related to the singularities of $X$.
	
	\begin{definition}
		The $d$-th Fitting ideal of $\Omega_X^1$ is a quasicoherent ideal in $X$ known as the jacobian ideal of $X$. The closed subscheme it defines will be called the \textbf{singular locus subscheme} of $X$ and denoted by $\Xs$. 
	\end{definition}
	Indeed, the underlying closed subset of $\Xs$ is the set of singular points of $X$.

	\begin{definition}
		Consider the natural $\O_Y$-morphism $h^\ast \Omega_X^d \to \Omega_Y^d$. Its image can be written as $\ \J_h \cdot \Omega_Y^d $ for a unique ideal $\J_h$ of $Y$. This is a coherent ideal, known as the \textbf{jacobian ideal} of $h$.
	\end{definition}

	\begin{remark} \normalfont
		Equivalently, $\J_h$ can be defined as the zeroth Fitting ideal of $\Omega_{Y/X}^1$.
	\end{remark}

	As explained in \cite[p. 522]{GW}, the underlying subset of $V(\J_h)$ is the support of $\Omega^1_{Y/X}$. Hence, it is contained in $h^{-1}(Z)$.
	
	We will use the following criterion to show that the subsets $\X_{m,E}$ are cylinders.
	
	\begin{proposition}[{\cite[Corollary 6.4]{EM}}] \label{criteria}
		Let $e,e'$ be two nonnegative integers and consider the subset
		\[
		C_{e,e'} \coloneqq \X_e(X,V(\J_h)) \ \cap \ h_\infty^{-1}\X_{e'}(X,\Xs)
		\]
		of $\L(Y)$. If $C$ is a cylinder of $\L(Y)$ contained in $C_{e,e'}$, then $h_\infty(C)$ is a cylinder of $\L(X)$.
	\end{proposition}

	\begin{proposition} \label{cylinder}
		If $h$ is $m$-separating, then for every $m$-divisor $E$, the subset $\X_{m,E}$ is an irreducible cylinder of $\L(X)$. Its closure in $\L(X)$ is a quasicylinder (see \cite[Definition 3.2]{FEI}) of codimension
		\[
		\frac{m(\widehat{k}_E+1)}{N_E}.
		\]
		Here, $\widehat{k}_E$ is the Mather discrepancy of $E$, see \cite[Definition 1.9]{FEI}.
	\end{proposition}
	\begin{proof}
		Consider the subset
		\[
		\widetilde{\X}_{m,E} \coloneqq \left\{\gt \in \L(Y) \ \middle| \ \gt \cdot E = \frac{m}{N_E}, \ \gt(0) \in E^\circ\right\},
		\]
		so $\X_{m,E} = h_\infty(\widetilde{\X}_{m,E})$. The fact that $h^{-1}(Z)$ is a simple normal crossing divisor allows us to check locally that $\widetilde{\X}_{m,E}$ is irreducible. Thus, so is $\X_{m,E}$.
		
		Recall that $\Xs$ is contained in $Z$ set-theoretically, so there exists $n \geq 1$ such that $\Xs$ is a closed subscheme of the $n$-th thickening $Z_n$ of $Z$. Hence, for every $\gt \in \widetilde{\X}_{m,E}$,
		\[
		\gt \cdot h^{-1}(\Xs) = h_\infty(\gt) \cdot \Xs \leq h_\infty(\gt) \cdot Z_n = (h_\infty(\gt) \cdot Z)n = mn,
		\]
		i.e. the intersection multiplicity of $\gt$ and $h^{-1}(\Xs)$ is bounded.
		
		Since $V(\J_h)$ is contained in $h^{-1}(Z)$ set-theoretically, a similar argument shows that the intersection multiplicity of any element of $\widetilde{\X}_{m,E}$ and $V(\J_h)$ is also bounded.
		
		For $e, e' \in \N$, denote
		\[
		\widetilde{\X}_{m,E;e,e'} \coloneqq \left\{\gt \in \widetilde{\X}_{m,E} \ \middle| \ \gt \cdot V(\J_h) = e, \ \gt \cdot h^{-1}(\Xs) = e' \right\},
		\]
		so
		\[
		\X_{m,E} = \bigcup_{e,e' \in \N} \widetilde{\X}_{m,E;e,e'}.
		\]
		By the discussion above, this union is finite. Also, note that the subset $\widetilde{\X}_{m,E;e,e'}$ is a locally closed cylinder in $\L(Y)$, so by Proposition \ref{criteria}, its image $h_\infty(\widetilde{\X}_{m,E;e,e'})$ is a cylinder of $\L(X)$. Hence the assertion follows.
		
		The codimension formula is deduced from combining Proposition \ref{generic-point} and \cite[Theorem 3.9]{FEI}.
	\end{proof}

	We have an immediate corollary, cf. \cite[Proposition 1.19]{BdlB} and \cite[Proposition 3.5]{Is13}.
	
	\begin{corollary}
		The codimension of the $m$-contact locus is given by
		\[
		\frac{1}{m}\codim\X_m(X,Z,\Sigma) = \inf \frac{\widehat{k}_E+1}{N_E},
		\]
		where the infimum runs over the $m$-divisors $E$ of $h$. In particular, the right-hand side does not depend on the particular choice of $m$-separating embedded resolution.
	\end{corollary}

	\subsection{The problem}
	Proposition \ref{cylinder} tells us that (\ref{union}) is a disjoint union of irreducible cylinders. In particular, every irreducible component of $\X_m$ is of the form $\overline{\X_{m,E}}$ for a unique $m$-divisor $E$ (here, bars denote the closure in $\X_m$). \\
	
	\noindent
	\textbf{embedded Nash problem.} Characterize the $m$-divisors $E$ such that $\overline{\X_{m,E}}$ is an irreducible component of $\X_m$. \\
	
	For the resolution of the problem, we will be interested in determining whether for a given pair of $m$-divisors $E$ and $E'$, the subset $\X_{m,E}$ is contained in $\overline{\X_{m,E'}}$ or not. Thanks to Proposition \ref{generic-point}, this can be restated in terms of maximal divisorial sets: $\X_{m,E}$ is contained in $\overline{\X_{m,E'}}$ if and only if $C_X(w_E)$ is contained in $C_X(w_{E'})$.

	\subsection{A topological necessary condition for adjacencies to happen}

	In the same philosophy as \cite{FdB} and \cite{BPP}, we give the following necessary condition for an adjacency to happen.
	
	\begin{lemma} \label{convergent-wedge}
		Let $E$ and $E'$ be two different $m$-divisors. If $\X_{m,E}  \subseteq \overline{\X_{m,E'}}$, then there exists a convergent $\C$-wedge $\alpha: \Spec \C[ \! [s] \! ] \to \L(X)$ such that $\alpha(0) \in \X_{m,E}$ and $\alpha(\eta) \in \X_{m,E''}$ for some $m$-divisor $E''$ such that $\X_{m,E''}  \not \subseteq \overline{\X_{m,E}}$.
	\end{lemma}
	
	This section is dedicated to proving the lemma above. To this end, we recall some definitions and results.
	
	\begin{definition}[{\cite[Definition 3.1]{Re}}]
		A closed irreducible subset $N$ of $\L(X)$ is called \textbf{generically stable} if it is not contained in $\L(\Xs)$ and if there exists an open affine subset $U$ such that $U \cap N$ is a nonempty closed subset of $U$ whose defining radical ideal is the radical of a finitely generated ideal.
	\end{definition}

	\begin{remark} \normalfont
		Strictly speaking, the above definition does not appear in \cite{Re}, but in the corrigendum \cite{Re21}. The results of \cite{Re} have been revised in \cite{Re21}.
	\end{remark}

	\begin{proposition} \label{gen-stable}
		Maximal divisorial sets are generically stable.
	\end{proposition}
	\begin{proof}
		Let $\mathfrak{Y}$ be an integral normal scheme, $\mathfrak{h}: \mathfrak{Y} \to X$ a proper birational morphism and $\mathfrak{E}$ a prime divisor in $\mathfrak{Y}$. Denote by $\val_\mathfrak{E}$ the divisorial valuation on $X$ associated to $\mathfrak{E}$. Let $q$ be a positive integer and set $v \coloneqq q\val_\mathfrak{E}$.
		
		By \cite[Proposition 3.4]{Is08}, $C_X(v)$ equals the closure of $\mathfrak{h}_\infty(\X_q(\mathfrak{Y},\mathfrak{E}))$ in $\L(X)$. Let $\gt$ be the generic point of the closure of $\X_q(\mathfrak{Y},\mathfrak{E})$, so $\gamma \coloneqq \mathfrak{h}_\infty(\gt)$ is the generic point of $C_X(v)$. We claim that $\gt$ is a fat schematic arc of $\mathfrak{Y}$. Indeed, note that $\gt(0)$ is the generic point of $\mathfrak{E}$, so $\mathfrak{E}$ is contained in the closure of $\gt(\eta)$. On the other hand, the fact that $\gt\cdot \mathfrak{E} = q < \infty$ means that $\gt(\eta) \not \in \mathfrak{E}$. Since $\mathfrak{E}$ is of codimension 1, the claim follows.
		
		Since $\mathfrak{h}$ is dominant, $\gamma$ is also a fat schematic arc of $X$. Now, the proposition follows from \cite[Lemma 3.6]{Re}.
	\end{proof}

	In \cite[Corollary 4.8]{Re}, Reguera gave a curve selection lemma for generically stable subsets in arc spaces. The proof was based on Corollary 4.6 in loc. cit., which claimed that the completion of the local ring of the arc space in a generically stable subset is noetherian. Later, the same result was proven by de Fernex and Docampo by analyzing the Kähler differentials of the arc space, see \cite[Corollary 10.13]{FD20}.
	
	We state the following version of the curve selection lemma in arc spaces.
	
	\begin{proposition}[{\cite[Lemma 6]{FdB}}] \label{csl}
		Let $N \subset N'$ be two different closed irreducible subsets of $\L(X)$. Assume that $N$ is generically stable. Let $\mathfrak{Z}$ be another closed subset in $\L(X)$ not containing $N'$. Let $\gamma$ be the generic point of $N$ and denote by $k$ its residue field. There exists a finite field extension $K/k$ and a $K$-wedge $\alpha$ in $X$ whose base point is $\gamma$ and whose generic point is in $N'-\mathfrak{Z}$.
	\end{proposition}
	
	We now discuss closed and $\C$-valued points in arc spaces. In a given $\C$-scheme, every $\C$-valued point is closed, and moreover, the converse also holds when the scheme is locally of finite type. Since our base field $\C$ is uncountable, we have an analogous result for schemes locally of countable type, which is due to Ishii. See also \cite[Ch. 3, Remark 3.3.11]{CNS}.
	
	\begin{proposition} \label{countable}
		In a $\C$-scheme of locally countable type, every closed point is a $\C$-valued point.
	\end{proposition}
	\begin{proof}
		Let $\mathfrak{X}$ be a $\C$-scheme of locally countable type and $x$ a closed point of $\mathfrak{X}$. Since the assertion is local, we may assume that $\mathfrak{X}$ is affine and a closed subscheme of $M \coloneqq \Spec \C[x_1,x_2,\dots]$. Note that $x$ is also a closed point of $M$, and the residue fields of $x$ in $\mathfrak{X}$ and in $M$ are isomorphic over $\C$. Thus, it suffices to show the claim for $M$. This is the content of the proof of \cite[Proposition 2.10]{Is04}. 
	\end{proof}

	We have the following consequence, cf. \cite[Proposition 3.35]{GW}.

	\begin{proposition} \label{very-dense}
		Every nonempty locally closed subset of a $\C$-scheme of locally countable type contains a $\C$-valued point.
	\end{proposition}
	\begin{proof}
		Let $\mathfrak{X}$ be a $\C$-scheme of locally countable type and $S$ a nonempty locally closed subset of $\mathfrak{X}$. Without loss of generality, we may assume that $S$ is a closed subset in an open affine subset $U$ of $\mathfrak{X}$. We write $U = \Spec A$ for a countable $\C$-algebra $A$ and $S = V(\mathfrak{a})$ for an ideal $\mathfrak{a}$ of $A$. The fact that $S$ is nonempty means that the ideal $\mathfrak{a}$ is proper, so it is contained in a maximal ideal of $A$. This maximal ideal defines a point $x$ in $S$ which is closed in $U$. Since $U$ is a $\C$-scheme locally of countable type, Proposition \ref{countable} implies that $x$ is a $\C$-valued point.
	\end{proof}
	
	We finally present a nested approximation theorem that generalizes Artin's classical result.
	
	\begin{proposition}[{\cite[Theorem 4.1]{BDLD}, \cite[Theorem 1.4]{Po}}] \label{approx}
		Consider algebraic power series
		\[
		f(x,y) = \big(f_1(x,y),\dots,f_N(x,y)\big) \in \C[ \! [x,y] \!]^N,
		\]
		where $x = (x_1,\dots,x_m)$ and $y = (y_1,\dots,y_n)$. Suppose that
		\[
		\overline{y}(x) = \big(\overline{y}_1(x),\dots,\overline{y}_n(x)\big) \in \C[ \! [x]\!]^n
		\]
		are formal power series such that $f\big(x,\overline{y}(x)\big) = 0$. Moreover, write $x' = (x_1,\dots,x_{m'})$ and $y' = (y_1,\dots,y_{n'})$ for $0 \leq m' \leq m$ and $0 \leq n' \leq n$, and assume that
		\[
		\big(\overline{y}_1(x),\dots,\overline{y}_{n'}(x)\big) \in \C[ \! [x']\!]^{n'}.
		\]
		Then, for every positive integer $c$, there exist algebraic power series
		\[
		y(x) = \big(y_1(x),\dots,y_n(x)\big) \in \C[ \! [x]\!]^n
		\]
		with
		\[
		\big(y_1(x),\dots,y_{n'}(x)\big) \in \C[ \! [x']\!]^{n'}
		\]
		such that
		$f\big(x,y(x)\big) = 0$ and
		\[
		y(x) \equiv \overline{y}(x) \mod \langle x \rangle^c.
		\]
	\end{proposition}

	\begin{proof}[Proof of Lemma \ref{convergent-wedge}]
		We will present a proof which adapts the arguments of sections 3 and 4 of \cite{FdB} to our setting.
		
		By Proposition \ref{gen-stable}, maximal divisorial sets are generically stable, so we can apply Proposition \ref{csl} where the closed subset we avoid is $C_X(w_E) \cup \X_{\geq m+1}$. Thus there exists a wedge $\alpha \in \L^{(2)}(X)$ such that $\alpha(0)$ is the generic point of $C_X(w_E)$ and $\alpha(\eta) \in C_X(w_{E'}) - (C_X(w_E) \cup \X_{\geq m+1})$. By Proposition \ref{generic-point}, we have that $\alpha(0) \in \X_{m,E}$. On the other hand, $\alpha(\eta) \in \X_m$, so $\alpha(\eta) \in \X_{m,E''}$ for some $m$-divisor $E''$. Since $\alpha(\eta) \not \in C_X(w_E)$, necessarily $C_X(w_{E''}) \not\subseteq C_X(w_E)$, i.e. $\X_{m,E''} \not\subseteq \overline{\X_{m,E}}$.
		
		Consider the set of schematic wedges
		\[
		\Lambda \coloneqq \{\alpha \in \L^{(2)}(X) \ | \ \alpha(0) \in \X_{m,E} \text{ and } \alpha(\eta) \in \X_{m,E''}\}.
		\]
		By the discussion above, it is nonempty. Moreover, taking base points and generic points define continuous maps $\L^{(2)}(X) \to \L(X)$, so by Proposition \ref{cylinder}, $\Lambda$ is a finite union of locally closed subsets. Therefore, $\Lambda$ contains a nonempty locally closed subset. The space of wedges of $X$ is a $\C$-scheme of countable type, so by Proposition \ref{very-dense} there exists a $\C$-wedge $\beta$ such that $\beta(0) \in \X_{m,E}$ and $\beta(\eta) \in \X_{m,E''}$.

		The wedge $\beta$ factors through an open affine subscheme of $X$, so without loss of generality we may assume that $X$ is affine. More precisely, we will assume that $X$ is the closed subscheme $V(f_1,\dots,f_N)$ of the affine space $\A^d$. The $\C$-wedge $\beta$ can be written as
		\[
		\beta(s,t) = \beta(s)(t) = \left(\sum_{j,k=0}^\infty b_{1jk}s^jt^k,\dots,\sum_{j,k=0}^\infty b_{djk}s^jt^k\right),
		\]
		and satisfies $f_n(\beta(s,t)) = 0$ for every $n=1,\dots,N$.
		
		Suppose that $\X_{m,E}$ and $\X_{m,E''}$ are cylinders of level $\ell$. There exist a finite number of regular functions $G_1,\dots,G_u,G',H_1,\dots,H_v,H'$ on $\L_\ell(X)$ such that
		\begin{align*}
			\beta(0) & \in \tau_\ell^{-1}\big(V(G_1,\dots,G_u)-V(G') \big) \subseteq \X_{m,E}, \\
			\beta(\eta) & \in \tau_\ell^{-1}\big(V(H_1,\dots,H_v)-V(H') \big) \subseteq \X_{m,E''}.
		\end{align*}
		Consider the field $\C(B_{ijk})$ of complex rational functions on the variables $\{B_{ijk}\}_{i=1,\dots,d, \ j \geq 0, \ k \geq 0}$, 
		and the $\C(B_{ijk})$-wedge
		\[
		\mathbb{B}(s,t) = \left(\sum_{j,k=0}^\infty B_{1jk}s^jt^k,\dots,\sum_{j,k=0}^\infty B_{djk}s^jt^k\right)
		\]
		in $X$. This should be thought of as a universal wedge in $X$.
		
		By assumption, $H'\big(\tau_\ell(\beta(\eta))\big)$ is a nonzero power series in $\Cs$. Denote by $\mu$ its $s$-adic order. In turn, the coefficient of the monomial $s^\mu$ in $H'\big(\tau_\ell(\mathbb{B}(\eta))\big)$ is a polynomial in $\C[B_{ijk} \ | \ i=1,\dots,d, \ j \geq 0, \ k = 0,\dots\ell]$. Let $M$ be an integer greater than $j$ if $B_{ijk}$ appears in such polynomial.
		
		Consider the system of equations
		\begin{align} \label{system}
			\begin{split}
				f_n\left(\sum_{k=0}^\ell B_{1k}t^k+t^{\ell+1}C_1,\dots,\sum_{k=0}^\ell B_{dk}t^k+t^{\ell+1}C_d\right) & = 0, \quad n=1,\dots,N \\
				H_q\left(\left(B_{ik}\right)_{\substack{i=1,\dots,d, \ k=0,\dots,\ell}}\right) & = 0, \quad q = 1,\dots,v
			\end{split} 
		\end{align}
		where the variables are $\{B_{ik}\}_{\substack{i=1,\dots,d, \ k=0,\dots,\ell}}$ and $\{C_i\}_{i=1,\dots,d}$, and take values in $\Cst$. Writing the wedge $\beta$ as
		\[
		\beta(s,t) = \left(\sum_{k=0}^\ell b_{1k}(s)t^k+t^{\ell+1}c_1(s,t),\dots,\sum_{k=0}^\ell b_{dk}(s)t^k+t^{\ell+1}c_d(s,t)\right),
		\]
		$B_{ik} = b_{ik}(s)$ and $C_i = c_i(s,t)$ is a solution of the system (\ref{system}). Applying Proposition \ref{approx}, we get algebraic power series $b'_{ik}(s) \in \Cs$ and $c'_i(s,t) \in \Cst$ such that $B_{ik} = b'_{ik}(s)$ and $C_i = c'_i(s,t)$ are solutions of (\ref{system}), and
		\begin{align*}
			b_{ik}(s) & \equiv b'_{ik}(s) \mod \langle s \rangle^{M+1}, \\
			c_i(s,t) & \equiv c'_i(s,t) \mod \langle s,t \rangle^{M+1}.
		\end{align*}
		We construct the $\C$-wedge $\beta'$ in $\A^d$ as
		\[
		\beta'(s,t) = \left(\sum_{k=0}^\ell b'_{1k}(s)t^k+t^{\ell+1}c'_1(s,t),\dots,\sum_{k=0}^\ell b'_{dk}(s)t^k+t^{\ell+1}c'_d(s,t)\right).
		\]
		It is an algebraic wedge in $X$ by construction. In particular, it is a convergent wedge. Also,
		\begin{align*}
			\tau_\ell(\beta'(0)) & = \left(\sum_{k=0}^\ell b'_{1k}(0)t^k,\dots,\sum_{k=0}^\ell b'_{dk}(0)t^k\right) = \left(\sum_{k=0}^\ell b_{1k}(0)t^k,\dots,\sum_{k=0}^\ell b_{dk}(0)t^k\right) = \\
			& = \tau_\ell(\beta(0)) \in V(G_1,\dots,G_u)-V(G'),
		\end{align*}
		so $\beta'(0) \in \X_{m,E}$.
		Moreover,
		\[
		\tau_\ell(\beta'(\eta)) = \left(\sum_{k=0}^\ell b'_{1k}(s)t^k,\dots,\sum_{k=0}^\ell b'_{dk}(s)t^k\right).
		\]
		It is in $V(H_1,\dots,H_v)$ by construction. Also, the coefficient of the monomial $s^\mu$ in $H'(\tau_\ell(\beta'(\eta)))$ coincides with the coefficient of the monomial $s^\mu$ in $H'(\tau_\ell(\beta(\eta)))$, which we know it is nonzero. In particular, $H'(\tau_\ell(\beta'(\eta))) \neq 0$, so altogether,
		\[
		\tau_\ell(\beta'(\eta)) \in V(H_1,\dots,H_v)-V(H'),
		\]
		and hence $\beta'(\eta) \in \X_{m,E''}$.
	\end{proof}

	\begin{remark} \normalfont
		The proof also shows that if the question posed in \cite[p. 127]{Re} has an affirmative answer, then in Lemma \ref{convergent-wedge}, the divisor $E''$ can be taken to be $E'$.
	\end{remark}

	\section{The case of surfaces}
	
	We treat the case when $X$ is a 
	surface, $Z = C$ is an effective Cartier divisor of $X$, and $\Sigma = \{o\}$, where $o \in C$ is a closed point. As in the previous sections, $h: Y \to X$ is assumed to be an $m$-separating embedded resolution of $(X,C,o)$. 
	As contact loci are local, we may assume that $X$ is affine and $C$ is principal if necessary.
	
	Given a divisor $E \in \CC$, its \textbf{valence} is defined to be the number of points in $E-E^\circ$, and will be denoted by $v_E$. A divisor $E \in \CC$ of valence one will be called a \textbf{leaf}.
	
	For the sake of notation, we will simply write $\CC_o$ for $\CC_{\{o\}}$. Note that every $E \in \CC_o$ is a smooth proper curve in $Y$, so its genus can be considered. We will denote it by $g_E$.
	
	We consider the following subsets of $\CC$. Let $L$ be a leaf. We denote by $\CC_L$ the set of divisors $E \in \CC$ such that there exist $n \geq 0$ and $E_0,E_1,\dots,E_n  \in \CC_o$ such that:
	\begin{enumerate}
		\item $E_0 = L$ and $E_n=E$
		\item For every $1 \leq i \leq n-1$, $E_i$ has valence two. 
		\item For every $0 \leq i \leq n-1$, $E_i$ has genus zero. 
		\item For every $1 \leq i \leq n$, $E_{i-1}$ and $E_i$ intersect. 
	\end{enumerate}
	Note that if $L \not \in \CC_o$, then $\CC_L = \emptyset$. Otherwise, $L \in \CC_L$. For instance, if $L \in \CC_o$ has positive genus, then $\CC_L = \{L\}$.
	
	The set $\CC_L$ has a natural total order. Namely, if $E,E' \in \CC_L$, then $E \leq E'$ if and only if there exist integers $0 \leq n \leq n'$ and $E_0,\dots,E_n,\dots,E_{n'} \in \CC_o$ such that:
	\begin{enumerate}
		\item $E_0 = L$, $E_n=E$ and $E_{n'}=E'$.
		\item For every $1 \leq i \leq n'-1$, $E_i$ has valence two. 
		\item For every $0 \leq i \leq n'-1$, $E_i$ has genus zero.
		\item For every $1 \leq i \leq n'$, $E_{i-1}$ and $E_i$ intersect.
	\end{enumerate}
	
	If $L \in \CC_o$, we denote by $E_L$ the maximal element of $\CC_L$. Note that $E_L$ is the only element in $\CC_L$ with either valence greater than two or positive genus.
	
	\begin{remark} \normalfont
		Alternatively, the subsets $\CC_L$ can be easily viasualised in the \emph{dual graph} of $h^{-1}(C)$, see \cite[\S 2.2.2]{Ne}. More precisely, consider the maximal simple path in the dual graph which starts at the vertex corresponding to $L$ and only contains vertices of valence two and genus zero. Then $\CC_L$ is the set of divisors corresponding to vertices in the path, together with the divisor that corresponds to the vertex that meets the last vertex of the path without being part of it, namely, $E_L$.
	\end{remark}
	
	The divisors of $\CC_L$ have nice numerical properties, which are recorded in the following lemma:
	
	\begin{lemma} \label{num-cl}
		Let $L \in \CC$ be a leaf.
		\begin{enumerate}
			\item For every $E \in \CC_L$, we have that $N_L$ divides $N_E$.
			\item If $E,E' \in \CC_L$ are two different intersecting divisors, then $\gcd\{N_E,N_{E'}\}=N_L$.
		\end{enumerate}
	\end{lemma}
	\begin{proof}
		(1) It will follow from (2).
		
		(2) We label the divisors of $\CC_L$ as $L=E_0 < E_1 < \ldots < E_n$. We show that $\gcd\{N_{E_{i-1}},N_{E_i}\} = N_L$ by induction on $i$.
		
		We assume that $C$ is principal, so we have the intersection number $C \cdot L = 0$. Since $L$ only meets $E_1$, this implies that
		\[
		L^2N_L+N_{E_1} = 0.
		\]
		In particular, $N_L$ divides $N_{E_1}$, and
		\[
		\gcd\{N_{E_0},N_{E_1}\} = \gcd\{N_L,N_{E_1}\} = N_L.
		\]
		Now suppose that $i>1$ and $\gcd\{N_{E_{i-2}},N_{E_{i-1}}\} = N_L$. Arguing again with intersection numbers, we have that
		\[
		N_{E_{i-2}}+E_{i-1}^2N_{E_{i-1}}+N_{E_i} = 0.
		\]
		Hence
		\[
		\gcd\{N_{E_{i-1}},N_{E_i}\} = \gcd\{N_{E_{i-2}},N_{E_{i-1}}\}=N_L.
		\]
	\end{proof}
	
	The set of $m$-divisors of $\CC_L$ will be denoted by $\CC_{L,m}$. A leaf $L$ will be called an \textbf{$m$-leaf} if it is also an $m$-divisor. By Lemma \ref{num-cl}(1), $\CC_{L,m}$ is nonempty if and only if $L$ is an $m$-leaf. In that case, we will denote the maximal element of $\CC_{L,m}$ by $E_{L,m}$. We remark that different $m$-leaves may have the same maximal element.
	
	Finally, need to impose some minimality condition on the embedded resolution.
	
	\begin{definition}
		We say that $h: Y \to X$ is \textbf{admissible} if for every divisor $E \in \CC$ of valence greater than two, there at at least three divisors $E' \in \NN_E$ such that $N_E$ does not divide $N_{E'}$.
	\end{definition}
	
	The following is the main result of the paper.
	
	\begin{theorem} \label{main}
		Let $X$ be a surface, $C$ an effective Cartier divisor, and $o$ a closed point in $C$. Suppose that $h$ is an admissible $m$-separating embedded resolution of $(X,C,o)$. The irreducible components of $\X_m(X,C,o)$ are the following:
		\begin{enumerate}
			\item $\overline{\X_{m,E_{L,m}}}$ for every $m$-leaf $L$.
			\item $\overline{\X_{m,E}}$ if $E$ is an $m$-divisor which does not lie in $\CC_{L,m}$ for any $m$-leaf $L$.
		\end{enumerate}
	\end{theorem}
	
	Theorem \ref{main} is a direct consequence of Theorem \ref{allowed} and Theorem \ref{forbidden} in the next sections. Namely, the former result will discard some candidate $m$-divisors, whereas the latter will show that the remaining ones give irreducible components.
	
	\subsection{Allowed adjacencies} \label{S-allowed}
	
	The first half of Theorem \ref{main} is the following result.
	
	\begin{theorem} \label{allowed}
		Let $L$ be an $m$-leaf, and let $E,E' \in \CC_{L,m}$. If $E \leq E'$, then $\X_{m,E} \subseteq \overline{\X_{m,E'}}$.
	\end{theorem}

	\begin{remark} \normalfont
		For Theorem \ref{allowed}, we do not need $h$ to be admissible.
	\end{remark}

	The proof will be done by reducing to the cyclic quotient surface singularity case. To do so, we will need the following contraction result, which is essentially an application of Artin's contraction theorem. 

	\begin{lemma} \label{artin-contr}
		Let $V$ be a smooth proper surface over $\C$ and $C$ a reduced connected curve in $V$, whose irreducible components are labelled as $C_1,\dots,C_n$. Suppose that the following conditions hold:
		\begin{enumerate}
			\item The intersection matrix $(C_i \cdot C_j)_{i,j}$ is negative-definite.
			\item The arithmetic genus of each $C_i$ is zero.
			\item The curves $\{C_i\}$ form a string, i.e. the only curves that $C_i$ intersects are $C_{i-1}$ and $C_{i+1}$.
		\end{enumerate}
		Then there exists a surface $V'$, a closed point $o' \in V'$ and a $\C$-morphism $\pi: V \to V'$ such that $\pi^{-1}(o')_{\mathrm{red}} = C$ and $\pi$ is an isomorphism over $V'-\{o'\}$.
	\end{lemma}
	\begin{proof}
		We denote the self-intersection of $C_i$ by $-e_i$. Applying Castelnuovo's contraction theorem repeatedly, we may assume that $e_i > 1$ for every $i$.
		
		We wish to apply \cite[Theorem 2.3]{Ar62}. It only remains to verify that the arithmetic genus of an effective nonzero divisor $Z= \sum_{i=1}^nr_iC_i$ is nonpositive. We proceed by induction on $\sum r_i$.
		
		If $\sum r_i = 1$, then $Z = C_i$ for some $i$, which has arithmetic genus zero by assumption.
		
		Now assume that $\sum r_i > 1$. Denote $I_0 \coloneqq \{i \ | \ r_i > 0\}$ and
		\[
		Z' \coloneqq \sum_{i \in I_0} C_i. 
		\]
		For each effective divisor $D$ in $V$ we will denote by $p_a(D)$ its arithmetic genus. Let $K$ be a canonical divisor of $V$. By the adjunction formula,
		\[
		p_a(D) = 1+\frac{1}{2}\left(D \cdot (D+K)\right).
		\]
		Since the arithmetic genus of each $C_i$ is zero, this implies that $E_i \cdot K = e_i-2$ for every $i$. Hence one obtains
		\[
		p_a(Z') = 1+|I_0|-|I_0'|,
		\]
		where $I_0' \coloneqq \{i \in \{1,\dots,n-1\} \ | \ r_i = r_{i+1} = 1\}$. Hence $p_a(Z') \leq 0$.
		
		In general, we have that $Z \geq Z'$. The case of $Z = Z'$ has just been settled. If $Z \ne Z'$, then $Z-Z'$ is an effective divisor and the arithmetic genus of $Z$ satisfies
		\[
		p_a(Z) = p_a(Z')+p_a(Z-Z')+Z' \cdot (Z-Z') -1.
		\]
		By induction hypothesis, $p_a(Z-Z')$ is nonpositive. Therefore, it is enough to show that $Z' \cdot (Z-Z')$ is also nonpositive. Indeed, for every $i$,
		\[
		Z' \cdot E_i \leq 2-e_i \leq 0,
		\]
		so $Z' \cdot (Z-Z') \leq  0$.
	\end{proof}

	\begin{proof}[Proof of Theorem \ref{allowed}]
		We label the elements of $\CC_L$ as $E_1,\dots,E_r$ so that $E_1 \leq \ldots \leq E_r$. Also denote by $-e_i$ the self-intersection of $E_i$ for $i=1,\dots,r$. We need to show that for every $i \leq j$, the subset $C_X(w_{E_i})$ is contained $C_X(w_{E_j})$.
		
		
		Without loss of generality, we assume that $C$ is principal. By Nagata's compactification theorem, there exist integral proper $\C$-schemes $\overline{X}$ and $\overline{Y}$ that respectively contain $X$ and $Y$ as open dense subschemes and a $\C$-morphism $\overline{h}: \overline{Y} \to \overline{X}$ extending $h$. Note that $\overline{Y}$ is a complete surface that is smooth at every point of the curve $\sum_{i=1}^{r-1} E_i$. Also, the curve is connected and its components have genus zero. Finally, since the intersection matrix of $h^{-1}(o)$ is negative-definite, so is the intersection matrix of $\sum_{i=1}^{r-1} E_i$. By Lemma \ref{artin-contr}, there exists a surface $X'$, a point $o' \in X'$ and a morphism $h': \overline{Y} \to X'$ such that $h'^{-1}(o')_{\mathrm{red}} = \sum_{i=1}^{r-1} E_i$ and $h'$ is an isomorphism over $X'-\{o'\}$. By \cite[Theorem 2.3.1]{Ne}, $(X',o')$ is a cyclic quotient surface singularity. That is, there exist an open affine neighborhood $U$ of $o'$ in $X'$ and coprime integers $0 < q < n$ such that $U$ is isomorphic to $\A^2/\mu_n$, where the action of $\mu_n \coloneqq \{\zeta \in \C \ | \ \zeta^n = 1\}$ on $\A^2$ is given by $\zeta \cdot (x,y) = (\zeta x, \zeta^qy)$.
		
		The integers $n$ and $q$ can be explicitly computed by writing the negative continued fraction
		\[
		[e_{r-1},\dots,e_1] = e_{r-1}-\cfrac{1}{e_{r-2}-\cfrac{1}{\ddots-\cfrac{1}{e_1}}}
		\]
		as $\tfrac{n}{q}$ in its irreducible form. Additionally, for every $i=2,\dots,r$ we introduce the positive coprime integers $n_i$, $q_i$ such that
		\[
		\frac{n_i}{q_i} = [e_{i-1},\dots,e_1]
		\]
		in its irreducible terms. Note that $n_{r-1}=q$ and $n_r = n$. We also set $n_0=0$ and $n_1=1$, so
		\begin{equation} \label{Formula}
		n_{i-1}-e_in_i+n_{i+1} = 0
		\end{equation}
		for every $i=1,\dots,r-1$.
		
		Since $C$ is principal, we have that $h^{-1}(C) \cdot E_i = 0$ for every $i$. Hence, setting $N_{E_0}=0$, it follows that
		\begin{equation} \label{formula}
			N_{E_{i-1}}-e_iN_{E_i}+N_{E_{i+1}} = 0
		\end{equation}
		for every $i=1,\dots,r-1$. Comparing (\ref{Formula}) and (\ref{formula}) we deduce that $N_{E_i} = N_{E_1}n_i$ for every $i=0,\dots,r$. In particular, $N_{E_{r-1}}=N_{E_1}q$ and $N_{E_r}=N_{E_1}n$.
		
		The open subvariety $U$ is isomorphic to the affine toric variety $\Spec \C[\sigma^\vee \cap \Z^2]$, where $\sigma$ is the cone in $\R^2$ generated by $(1,0)$ and $(q,n)$. Moreover, $h'^{-1}(U)$ is a toric variety and $h': h'^{-1}(U) \to U$ is a toric morphism. Furthermore, there exists a finite toric morphism $U \to \A^2$, which ramifies over the axes and pullbacks one of the axes to $E_r$. Therefore, for every $i=1,\dots,r$ the valuation $w_{E_i}$ is a toric valuation on $U$.
		
		Therefore, we can apply \cite[Lemma 3.11]{Is08}; that is, the inclusion $C_{U}(w_{E_i}) \subseteq C_{U}(w_{E_j})$ holds if and only if
		\begin{equation} \label{ineq}
			w_{E_i}(f) \geq w_{E_j}(f)
		\end{equation}
		for every $f \in \O_{X'}(U) = \C[x,y]^{\mu_n}$. Once this is checked, \cite[Proposition 2.9(ii)]{Is08} implies the inclusion $C_{X'}(w_{E_i}) \subseteq C_{X'}(w_{E_j})$. 
		
		Furthermore, we claim that the rational map $\overline{h} \circ h'^{-1}: X' \dashrightarrow \overline{X}$ is in fact a morphism. Indeed, it is defined at every point different from $o'$ because $h'$ is an isomorphism over $X'-\{o'\}$. Since $\{o'\}$ has codimension 2 in $X'$ and $X'$ is normal, $\overline{h} \circ h'^{-1}$ extends uniquely to a morphism on $X'$. Since $X'$ is proper over $\C$, the morphism $\overline{h} \circ h'^{-1}$ is automatically proper. Therefore, if $C_{X'}(w_{E_i}) \subseteq C_{X'}(w_{E_j})$ holds, then by \cite[Proposition 2.9(i)]{Is08} so does the inclusion $C_X(w_{E_i}) \subseteq C_X(w_{E_j})$.
		
		Thus, it only remains to show inequality (\ref{ineq}). In \cite[\S 2.3.5]{Ne}, a minimal set of generators of $\C[x,y]^{\mu_n}$ is given. In summary, consider the convex hull in $\R^2$ of $(\sigma \cap \Z^2)-\{0\}$. Let $\{(p_j,q_j)\}_{j=1}^s$ be the sequence of lattice points in the boundary starting in $(p_0,q_0)=(1,0)$ and ending in $(p_s,q_s)=(q,n)$.
		
		Fix $j=1,\dots,s$. We set $v_j^r \coloneqq q_j$ and $v_j^{r-1} \coloneqq p_j$, and inductively,
		\[
		v_j^{i-1}-e_iv_j^i+v_j^{i+1} = 0
		\]	
		for $i=2,\dots,r-1$. We set $f_j(x,y) \coloneqq x^{v^0_j}y^{q_j}$. Then $\{f_j\}_{j=1}^s$ is a minimal set of generators of $\C[x,y]^{\mu_n}$. Moreover, in loc. cit. it is shown that $v_{E_i}(f_j) = v^i_j$.
		
		Since the $v_{E_i}$ are toric valuations, inequality (\ref{ineq}) only has to be checked for the generators $\{f_j\}_{j=1}^s$. It will follow from the following claim:
		\begin{proof}[Claim] \let\qed\relax
			For every $i=1,\dots,r-1$ and every $j=1,\dots,s$, we have that
			\begin{equation} \label{final-ineq}
			\frac{v^i_j}{N_{E_i}} \geq \frac{v^{i+1}_j}{N_{E_{i+1}}}.
			\end{equation}
		\end{proof}
		\begin{proof}[Proof of claim]  \let\qed\relax
			 We proceed by induction on $i$. If $i=r-1$, then the inequality becomes
			 \[
			 \frac{p_j}{N_{E_1}q} \geq \frac{q_j}{N_{E_1}n},
			 \]
			 i.e. $np_j \geq qq_j$. By construction, the point $(p_j,q_j)$ lies in the region of $\R^2$ defined by $nx \geq qy$, so the desired inequality holds.
		 	
		 	We now assume the inequality is true for $i+1$, $1 \leq i \leq r-2$. Writing $v^i_j = e_{i+1}v^{i+1}_j-v^{i+2}_j$ and $N_{E_i} = e_{i+1}N_{E_{i+1}}-N_{E_{i+2}}$, one sees that inequality (\ref{final-ineq}) is equivalent to
		 	\[
		 	\frac{v^{i+1}_j}{N_{E_{i+1}}} \geq \frac{v^{i+2}_j}{N_{E_{i+2}}}.
		 	\]
		 	In turn, this holds by induction hypothesis.
		\end{proof}
	\end{proof}

	\subsection{Forbidden adjacencies} \label{S-forbidden}
	
	The second half of Theorem \ref{main} is the following result.
	
	\begin{theorem} \label{forbidden}
		Suppose that $h$ is admissible, and let $E$ and $E'$ be two different $m$-divisors. If $E = E_{L,m}$ for some $m$-leaf $L$, or if $E \not \in \CC_{L,m}$ for every $m$-leaf $L$, then $\X_{m,E} \not\subseteq \overline{\X_{m,E'}}$.
	\end{theorem}
	
	The proof will be similar to the one in \cite[Theorem 7.4]{BdlB}. For us, Lemma \ref{convergent-wedge} will be the primary tool. Since the argument is topological, we recall the main construction used in it.
	
	\subsubsection{The A'Campo space} \label{acampo} Assume that $X$ is affine and $C$ is a principal divisor. Thus, $X \subseteq \A^N$ and there exists $f \in \O(X)$ such that $C = \mathrm{div} f$. By the Lê-Milnor theorem \cite[Theorem 1.1]{Le}, there exist $0 < \delta \ll \varepsilon \ll 1$ such that $f: B_\varepsilon \cap f^{-1}(S^1_\delta) \to S^1_\delta$ is a smooth locally trivial fibration. Let $\ft \coloneqq f \circ h \in \O(Y)$. Since $h: Y-h^{-1}(C) \to X-C$ is an isomorphism, this fibration is isomorphic to $\ft: h^{-1}(B_\varepsilon \cap f^{-1}(S^1_\delta)) \to S^1_\delta$.
	
	Consider the analytic map $\ft: h^{-1}(B_\varepsilon \cap f^{-1}(D_\delta)) \to D_\delta$ and the A'Campo construction \cite[\S 2]{A'C} associated to it. This is a topological manifold with boundary $T_{\AC}$ together with a fibration $f_{\AC}: T_{\AC} \to S^1$ isomorphic to $\ft: h^{-1}(B_\varepsilon \cap f^{-1}(S^1_\delta)) \to S^1_\delta$, i.e. there is a homeomorphism $a: h^{-1}(B_\varepsilon \cap f^{-1}(S^1_\delta)) \to T_{\AC}$ preserving the fibrations over $S^1$. Moreover, $T_{\AC}$ comes with a continuous map $\pi: T_{\AC} \to h^{-1}(B_\varepsilon \cap C)$.
	
	Denote by $F_{\AC,\theta}$ the fiber of $f_{\AC}$ over $\theta \in S^1$, and set $F_{\AC} \coloneqq F_{\AC,1}$. Given $E \in \CC$, we denote
	\[
	F_{\AC,E} \coloneqq F_{\AC} \cap \pi^{-1}(E^\circ \cap h^{-1}(B_\varepsilon)).
	\]
	The restriction $\pi: F_{\AC,E} \to E^\circ$ is an unramified covering of degree $N_E$.
	
	Similarly, given two different intersecting divisors $E, E' \in \CC$, we denote
	\[
	F_{\AC,\{E,E'\}} \coloneqq F_{\AC} \cap \pi^{-1}(E \cap E' \cap h^{-1}(B_\varepsilon)).
	\]
	These subsets give rise to a set-theoretically disjoint union decomposition of $F_\AC$ as
	\[
	F_\AC = \bigsqcup_{E \in \CC} F_{\AC,E} \ \sqcup \bigsqcup_{E \cap E' \neq \emptyset} F_{\AC,\{E,E'\}}.
	\]
	
	\begin{remark} \normalfont \label{top-pieces}
		Although it will not be essential in the rest of the paper, more can be said about the topology of these pieces. Let $E \in \CC$, and denote by $c_E$ the number of connected components of $F_{\AC,E}$. By \cite{Ne00}, $c_E$ divides $\gcd\{N_{E'} \ | \ E' \in \NN_E \cup \{E\}\}$. Moreover, if $E \in \CC_o$ and has genus zero, then both numbers coincide.
		
		If $E \in \CC_o$, then a direct application of the Riemann-Hurwitz formula shows that each connected component of $\overline{F_{\AC,E}}$ is a compact orientable surface of genus
		\[
		1+\frac{1}{c_E}\left(N_E\left(\frac{v_E}{2}+g_E-1\right)-\frac{1}{2}\left(\sum_{E' \in \NN_E}\gcd\{N_E,N_{E'}\}\right)\right)
		\]
		and
		\[
		\frac{1}{c_E}\sum_{E' \in \NN_E} \gcd\{N_E,N_{E'}\}
		\]
		boundary components.
		
		If $E \not \in \CC_o$, then $F_{\AC,E}$ has $v_E$ connected components, and each of them is a cylinder.
		
		Similarly, given two different intersecting divisors $E, E' \in \CC$, $F_{\AC,\{E,E'\}}$ is homeomorphic to a disjoint union of $\gcd\{N_E,N_{E'}\}$ cylinders.
	\end{remark}
	
	\begin{remark} \normalfont
		The considerations of Remark \ref{top-pieces} become particularly simple for divisors in $\CC_L$ for a leaf $L$. By Lemma \ref{num-cl}(2), $c_E=N_L$ for every $E \in \CC_L-\{E_L\}$. Moreover, if $E \neq L$, then $\overline{F_{\AC,E}}$ is homeomorphic to a disjoint union of $N_L$ cylinders. Also, if $L \neq E_L$, then $F_{\AC,E}$ is homeomorphic to a disjoint union of $N_L$ disks.
	
		Note that gluing a disk and a cylinder along the boundary of the disk and a boundary component of the cylinder produces a new disk, so altogether,
		\[
		\bigcup_{\substack{E \in \CC_L \\ E \neq E_L}} F_{\AC,E} \ \cup \bigcup_{\substack{E,E' \in \CC_L \\ E \cap E' \neq \emptyset}} F_{\AC,\{E,E'\}}
		\]
		has $N_L$ connected components, each of them homeomorphic to a disk.
	\end{remark}
	
	The fibration $f_\AC: T_\AC \to S^1$ comes with an explicit monodromy trivialization $\{\varphi_\tau: F_\AC \to F_{\AC,e^{2 \pi i \tau}}\}_{\tau \in \R}$. We set $\varphi \coloneqq \varphi_1 : F_\AC \to F_\AC$. The subsets of $F_{\AC,E}$ and $F_{\AC,\{E,E'\}}$ described above are invariant by $\varphi$. Moreover, the restriction $\varphi: F_{\AC,E} \to F_{\AC,E}$ is a covering transformation of $\pi: F_{\AC,E} \to E^\circ$, and in fact it generates the group of covering transformations of $\pi$. Also importantly, the set of fixed points of $\varphi^m$ is
	\[
	\operatorname{Fix} \varphi^m = \bigsqcup_{N_E | m} \overline{F_{\AC,E}}.
	\]

	\begin{lemma} \label{topological-lemma}
		Let $\eta: [0,1] \to F_\AC$ be a path whose endpoints are fixed by $\varphi^m$. If $\varphi^m \circ \eta$ and $\eta$ are homotopic relative to the endpoints, then either:
		\begin{itemize}
			\item There exists an $m$-divisor $E$ such that both endpoints of $\eta$ are in $\overline{F_{\AC,E}}$.
			\item There exists an $m$-leaf $L$ and two $m$-divisors $E,E' \in \CC_L$ such that $\eta(0) \in \overline{F_{\AC,E}}$ and $\eta(1) \in \overline{F_{\AC,E'}}$.
		\end{itemize}
	\end{lemma}
	\begin{proof}
		For simplicity, for every $E \in \CC$, we will consider the manifolds $E^\circ$ and $F_{\AC,E}$ as compact manifolds with boundary, instead of noncompact manifolds with empty boundary.
		
		Let $E \in \CC$. For each $E' \in \NN_E$, pick a point $p_{E'}=p_{E,E'}$ in the boundary component of $E^\circ$ corresponding to $E'$. Also, pick a point $p=p_E$ in the interior of $E^\circ$. Let $\alpha_{E'}=\alpha_{E,E'}$ be the loop based at $p_{E'}$ around the boundary of $E^\circ$ corresponding to $E'$. Also, let $\beta_{E'}=\beta_{E,E'}$ be a path from $p$ to $p_{E'}$. Finally, let $\sigma_1,\tau_1,\dots,\sigma_g,\tau_g$ the usual generators of the fundamental group of $E$ based at $p$. The fundamental groupoid of $E^\circ$ with base points $\{p_{E'}\} \cup \{p\}$ is generated by $\{\alpha_{E'}\}_{E'} \cup \{\beta_{E'}\}_{E'} \cup \{\sigma_1,\tau_1,\dots,\sigma_g,\tau_g\}$, subject to the unique relation
		\begin{equation}\label{relation}
			\prod_{E' \in \NN_E} \beta_{E'} \alpha_{E'} \beta_{E'}^{-1} = \prod_{j=1}^g \sigma_j\tau_j\sigma_j^{-1}\tau_j^{-1}, 
		\end{equation}
		in an appropriate ordering of $\NN_E$.
		
		Pick a preimage $p_0=p_{E,0}$ of $p$ by the covering map $\pi: F_{\AC,E} \to E^\circ$. For each $i \in \Z/N_E$, define $p_i = p_{E,i} \coloneqq \varphi_{\AC}^i(p_0)$. Denote by $\beta_{E',i}=\beta_{E,E',i}$ the lift of $\beta_{E'}$ starting at $p_i$, and denote by $p_{E',i}=p_{E,E',i}$ its ending point. Note that $\varphi^i_{\AC} \circ \beta_{E',0} = \beta_{E',i}$, so in particular $\varphi^i_{\AC}(p_{E',0})=p_{E',i}$. The lift $\alpha_{E',i}=\alpha_{E,E',i}$ of $\alpha_{E'}$ starting at $p_{E',i}$ ends at $p_{E',i-N_{E'}}$. Similarly, denote by $\sigma_{j,i}$ and $\tau_{j,i}$ the lifts of $\sigma_j$ and $\tau_j$ starting at $p_i$, respectively. The fundamental groupoid of $F_{\AC,E}$ with base points $\{p_i\}_i \cup \{p_{E',i}\}_{E',i}$ is generated by $\{\alpha_{E',i}\}_{E',i} \cup \{\beta_{E',i}\}_{E',i} \cup \{\sigma_{j,i}, \tau_{j,i}\}_{i,j}$, subject to the $N_E$ relations coming from lifting (\ref{relation}).
		
		Suppose that $E$ has valence greater than two or positive genus, and let $E' \in \NN_E$. We claim that the loop
		\[
		\prod_{k=0}^{N_E/\gcd\{N_E,N_{E'}\}-1} \alpha_{E',-N_{E'}k}
		\]
		has infinite order in the fundamental grupoid of $F_{\AC,E}$ even after adding for each $E'' \in \NN_E-\{E'\}$ the relation
		\[
		\prod_{k=0}^{N_E/\gcd\{N_E,N_{E''}\}-1} \alpha_{E'',-N_{E''}k} = 1.
		\]
		Note that this is equivalent to showing that $\alpha_{E'}^{N_E/\gcd\{N_E,N_{E'}\}}$ has infinite order in the fundamental groupoid of $E^\circ$ after adding for each $E'' \in \NN_E-\{E'\}$ the relation
		\[
		\alpha_{E''}^{N_E/\gcd\{N_E,N_{E''}\}} = 1.
		\]
		After ordering $\NN_E$ appropriately, we can assume that
		\[
		\beta_{E'}\alpha_{E'}\beta_{E'}^{-1} = \prod_{E'' \in \NN_E-\{E'\}} \beta_{E''} \alpha_{E''} \beta_{E''}^{-1} \cdot \prod_{j=1}^g \sigma_j\tau_j\sigma_j^{-1}\tau_j^{-1}.
		\]
		If $E$ has positive genus, then the right-hand side has infinite order. Now, suppose that $E$ has genus zero. Since $h$ is admissible and $E$ has valence greater than two, there exist at least two divisors $E'' \in \NN_E-\{E'\}$ such that $N_E$ does not divide $N_{E''}$. This means that $N_E > \gcd\{N_E,N_{E''}\}$, and hence, the loops $\alpha_{E''}$ are nontrivial. This implies again that the right-hand side has infinite order. Therefore, in both cases, $\alpha_{E'}$ has infinite order, and in particular, so does $\alpha_{E'}^{N_E/\gcd\{N_E,N_{E'}\}}$.
		
		On the other hand, the fundamental groupoid of $F_{\AC,\{E,E'\}}$ with base points $\{p_{E,E',i}\}_i \cup \{p_{E',E,i}\}_i$ is free generated by $\{\alpha_{E,E',i}\}_i \cup \{\gamma_{E,E',i}\}_i$. Here, $\gamma_{E,E',i}$ is a path from $p_{E,E',i}$ to $p_{E',E,i}$.
		
		By the Seifert-Van Kampen theorem for groupoids \cite[\S 9.1]{Br}, we obtain the fundamental groupoid of $F_\AC$ over all the previous base points. Of course, one can see this groupoid is free by deleting extra generators of the form $\alpha_{E',i}$.
		
		Finally, let $\eta$ be a path between two fixed points of $\varphi^m$ which lie on different connected components of $\operatorname{Fix} \varphi^m$. Without loss of generality we may assume that the endpoints are some of the base points of the fundamental groupoid of $F_\AC$. We suppose that there is no leaf $L$ such that $\eta(0) \in \overline{F_{\AC,E}}$ and $\eta(1) \in \overline{F_{\AC,E'}}$ for some $m$-divisors $E,E' \in \CC_L$, so we need to show that $\eta$ and $\varphi^m \circ \eta$ are not homotopic.
		
		We write $\eta$ as a product on the generators. This expression can be taken to be unique if we choose a subset of generators over which the groupoid is free. We distinguish two cases:
		\begin{itemize}
			\item Suppose that there is a generator $\gamma_{E,E',i}$ in the expression of $\eta$ such that there is no leaf $L$ such that $E,E' \in \CC_L$. If $N=\lcm\{N_E,N_{E'}\}$ and
			\[
			c=	\prod_{k=0}^{N_E/\gcd\{N_E,N_{E'}\}-1} \alpha_{E,E',i-N_{E'}k},
			\]
			then $(\varphi^m)^N(\gamma_{E,E',i}) = c \cdot \gamma_{E,E',i}$. Since $c$ has infinite order, we conclude that $(\varphi^m)^N \circ \eta$ and $\eta$ are not homotopic. In particular, nor are $\varphi^m \circ \eta$ and $\eta$.
			\item If there is no such generator $\gamma_{E,E',i}$ in the expression of $\eta$, then there exists a generator $\beta_{E,E',i}$ in the expression of $\eta$ corresponding to a divisor $E$ which is not an $m$-divisor. 
			Therefore, $\varphi^m \circ \beta_{E',i} = \beta_{E',i+m} \neq \beta_{E',i}$. Hence $\varphi^m \circ \eta$ and $\eta$ are not homotopic.
		\end{itemize}
	\end{proof}
	
	\begin{lemma} \label{no-convergent-wedge}
		If $E$ and $E'$ are two different $m$-divisors such that there is no $m$-leaf $L$ such that $E,E' \in \CC_{L,m}$, then there is no convergent $\C$-wedge $\alpha: \Spec \C [\![s]\!] \to \L(X)$ such that $\alpha(0) \in \X_{m,E}$ and $\alpha(\eta) \in \X_{m,E'}$.
	\end{lemma}
	\begin{proof}
		We will reproduce the arguments of the proof of \cite[Theorem 7.4]{BdlB}.
		
		Without loss of generality, we assume that $X$ is affine and $C$ is principal. Arguing by contradiction, suppose that such a $\C$-wedge $\alpha$ exists. Let $f \in \mathcal{O}(X)$ such that $C = \mathrm{div} f$. There exists a convergent $\C$-wedge $u: \Spec \C [ \! [s] \! ] \to \L(\A^1)$ on $\A^1$ such that $f(\alpha(s))=u(s)^m$. Moreover, the power series $u$ is of $t$-adic order 1. Therefore, we can make the analytic change of coordinates $(\tilde{s},\tilde{t}) = (s,u(s)(t))$ so that $f(\alpha(\tilde{s},\tilde{t}(s,t))) = \tilde{t}^m$. Therefore, we can assume that $f(\alpha(s)) = t^m$.
		
		Suppose that $\alpha$ is defined on $[0,1]$. Thus, for every $s_0 \in [0,1]$ we may replace $s$ by $s_0$ in $\alpha$ and get a convergent $\C$-arc, which we denote by $\alpha(s_0)$. We set $E_0 \coloneqq E$ and $E_1 \coloneqq E'$ to shorten the notation. For $s \in \{0,1\}$ write $N_s \coloneqq N_{E_s}$ and $p_s \coloneqq \widetilde{\alpha(s)}(0) \in E_s$. Also, denote by $(U_s,\psi_s=(x_s,y_s))$ a complex analytic chart of $p_s$ contained in $h^{-1}(B_\varepsilon \cap f^{-1}(D_\delta))$ such that $\psi_s(p_s) = (0,0)$, and $f \circ \psi_s^{-1} = x_s^{N_s}$. Let $\widetilde{\alpha(s)} = (x_s(t),y_s(t))$ in the coordinates of $(U_s,\varphi_s)$. We have that $x_s(t)^{N_s} = t^m$, so without loss of generality $x_s(t) = t^{m/N_s}$. Consider the path $[0,1] \ni \lambda \mapsto (t^{m/N_s},\lambda y_s(t))$ in $\X_m$ between $\widetilde{\alpha(s)}$ and $(t^{m/N_s},0)$. Concatenating with $\alpha$, this procedure gives us a path $\beta: [0,1] \to \X_m$ such that $\beta(s)$ is convergent and $f(\beta(s))= t^m$ for all $s \in [0,1]$, and $\widetilde{\beta(0)} = (t^{m/N_E},0)$ and $\widetilde{\beta(1)} = (t^{m/N_{E'}},0)$.
		
		Consider A'Campo's construction described in \S \ref{acampo} with the same notation. As explained in the proof of \cite[Theorem 7.4]{BdlB}, the homeomorphism $a$ considered above can be chosen so that for $s \in \{0,1\}$,
		\begin{align*}
			a\big(\psi_s^{-1}(\delta^{1/N_s},0)\big) & \in F_{\AC,E_s}, \\
			\varphi_\tau\Big(a\big(\psi_s^{-1}(x_s,y_s)\big)\Big) & = a\big(\psi_s^{-1}(x_se^{2\pi i \tau/N_s},y_s)\big)
		\end{align*}
		Suppose that $\beta(s)$ is defined on $[0,\delta^{1/m}]$ for every $s \in [0,1]$. Consider the continuous map
		\[
		H: [0,1] \times [0,1] \to F_\AC; \quad (s,\tau) \mapsto \varphi_{m\tau}^{-1}\Big(a\big(\widetilde{\beta(s)}(\delta^{1/m}e^{2\pi i \tau})\big)\Big) \ .
		\]
		If $s \in \{0,1\}$, then
		\[
		H(s,\tau) = \varphi_{m\tau}^{-1}\Big(a\big(\widetilde{\beta(s)}(\delta^{1/m}e^{2\pi i \tau})\big)\Big) = a\big(\psi_s^{-1}(\delta^{1/N_s},0)\big) \in F_{\AC,E_s}.
		\]
		Therefore, $H$ is a homotopy between $H_0 \coloneqq H(\cdot,0)$ and $H_1 \coloneqq H(\cdot,1)$ relative to the endpoints. Moreover, $\varphi^m \circ H_1 = H_0$ by construction. However, $H_0$ contradicts Lemma \ref{topological-lemma}.
	\end{proof}
	
	\begin{proof}[Proof of Theorem \ref{forbidden}]
		Arguing by contradiction, suppose that $\X_{m,E} \subseteq \overline{\X_{m,E'}}$. By Lemma \ref{convergent-wedge}, there exists a convergent $\C$-wedge $\alpha: \Spec \C[ \! [s] \! ] \to \L(X)$ such that $\alpha(0) \in \X_{m,E}$ and $\alpha(\eta) \in \X_{m,E''}$ for some $m$-divisor $E''$ such that $\X_{m,E''}  \not \subseteq \overline{\X_{m,E}}$.
		
		If $E = E_{L,m}$, then Theorem \ref{allowed} in particular implies that $E'' \not \in \CC_{L,m}$. Therefore, the result follows from applying Lemma \ref{no-convergent-wedge} to $E$ and $E''$.
	\end{proof}

	\subsection{Comparison with the smooth-ambient case}
	
	Theorem \ref{main} recovers \cite[Theorem 1.22]{BdlB} when $X$ is assumed to be smooth. Moreover, in that case, the divisors $E \in \CC_o$ have all genus zero, so the description of the sets $\CC_L$ is more straightforward. Nevertheless, \cite[Theorem 1.22]{BdlB} also includes a second result: that the irreducible components of $\X_m$ are disjoint. Although the proof of Theorem \ref{main} does not clarify this point, 
	the evidence shows that this may still be true. We state it as a conjecture.
	
	\begin{conjecture} \label{conj}
		Let $X$ be an integral separated 
		surface, $C$ an effective Cartier divisor in $X$, and $o$ a closed point in $C$. Assume that $X$ is smooth away from $C$. Then the irreducible components of $\X_m(X,C,o)$ are disjoint.
	\end{conjecture}

	\vspace{3cm}
	
	\noindent
	\textsc{Javier de la Bodega}
	\vspace{0.2cm} \\
	Basque Center for Applied Mathematics. Alameda Mazarredo 14, 48009 Bilbao, Spain. \url{jdelabodega@bcamath.org}
	\vspace{0.2cm} \\
	KU Leuven, Department of Mathematics. Celestijnenlaan 200B, 3001 Heverlee, Belgium. \url{javier.delabodega@kuleuven.be}
		
\end{document}